\documentclass{amsart} \usepackage{amsmath} \usepackage{amscd}
\usepackage{amsxtra}

\usepackage{amssymb}
\usepackage[mathscr]{eucal}	
\DeclareSymbolFont{pssymbols}     {OMS}{ztmcm}{m}{n}
\DeclareSymbolFontAlphabet{\mathpsscr}   {pssymbols}

\usepackage{amsthm}
\theoremstyle{plain}
\newtheorem{thm}{Theorem}
\newtheorem{cor}[thm]{Corollary}
\newtheorem{prop}[thm]{Proposition}
\newtheorem{lem}[thm]{Lemma}
\newtheorem*{thm*}{Theorem}
\newtheorem*{cor*}{Corollary}
\newtheorem*{prop*}{Proposition}
\newtheorem*{lem*}{Lemma}

\theoremstyle{definition}
\newtheorem{defn}[thm]{Definition}
\newtheorem{examples}[thm]{Examples}

\theoremstyle{remark}
\newtheorem{rem}[thm]{Remark}
\newtheorem{rems}[thm]{Remarks}

\newtheorem*{rem*}{Remark}
\newtheorem*{rems*}{Remarks}
\newtheorem*{note*}{Note}

\newcommand{\itemref}[1]{\textup{(\ref{#1})}}
\usepackage[arrow,matrix,tips,frame,curve,ps,dvips]{xy}
\SelectTips{cm}{}
\UseTips
\CompileMatrices

\newdir^{ (}{{}*!/-5pt/@^{(}}
\newdir_{ (}{{}*!/-5pt/@_{(}}

\newcommand{\dash}{\nobreakdash-\hspace{0pt}}

\newcommand{\CC}{{\mathbb C}}
\newcommand{\RR}{{\mathbb R}}
\newcommand{\ZZ}{{\mathbb Z}}

\newcommand{\QQ}{{\mathbb Q}}

\DeclareMathOperator{\rank}{rank}

\renewcommand{\Re}{\operatorname{Re}}   
\newcommand{\PP}{\mathbb P}

\DeclareMathOperator{\QQrank}{\QQ-rank}
\DeclareMathOperator{\RRrank}{\RR-rank}
\DeclareMathOperator{\CCrank}{\CC-rank}
\DeclareMathOperator{\fieldrank}{-rank}

\DeclareMathOperator{\supp}{supp}

\renewcommand{\l}{\ell}

\DeclareMathOperator{\restr}{\rho}	
\newcommand{\rvvv}[1][]{\rVert_{\if!#1!b\else#1,b\fi}}	
\renewcommand{\setminus}{-}

\newcommand{\lsb}[1]{{}_{#1}}
\newcommand{\lsp}[1]{{}^{#1}\!}

\newbox\arrowbox
\setbox\arrowbox=\hbox{\lower .5ex\rlap{$\longrightarrow$}%
 \raise .5ex\hbox{$\longleftarrow$}}

\newcommand{\G}{\Gamma}

\DeclareMathOperator{\GL}{GL}
\DeclareMathOperator{\PGL}{PGL}

\DeclareMathOperator{\Gal}{Gal}

\newcommand{\Xhat}{\widehat{X}}
\newcommand{\back}{\backslash}
\newcommand{\Dstar}{D^*}
\newcommand{\Xstar}{X^*}

\newcommand{\al}{\alpha}

\renewcommand{\b}{\beta}
\renewcommand{\d}{\delta}

\newcommand{\D}{\Delta}

\newcommand{\DCC}{\lsb\CC\D}

\newcommand{\e}{\varepsilon}

\newcommand{\z}{\zeta}

\newcommand{\U}{\varUpsilon}

\newcommand{\kap}{\kappa}
\renewcommand{\o}{\omega}

\makeatletter
\def\sphat{^{\mathchoice{}{}%
 {\,\,\smash[b]{\hbox{\lower4\ex@\hbox{$\m@th\widehat{\null}$}}}}%
 {\,\smash[b]{\hbox{\lower3\ex@\hbox{$\m@th\hat{\null}$}}}}}\,}
\makeatother

\newcommand{\F}{\mathscr F}
\newcommand{\Ftilde}{\widetilde{\F}}
\newcommand{\B}{\mathscr B}
\newcommand{\K}{\mathscr K}

\begin{document}


\author{Leslie Saper}
\address{Department of Mathematics\\ Duke University\\ Box 90320\\ Durham,
NC 27708\\U.S.A.}
\email{saper@math.duke.edu}
\urladdr{http://www.math.duke.edu/faculty/saper}
\title[Geometric rationality]{Geometric rationality of equal-rank Satake compactifications}
\dedicatory{Dedicated to the memory of Armand Borel, with respect and
  admiration}
\thanks{
This research was supported in part by the  National Science Foundation
through  grant  DMS-9870162 and through support to attend the 2002 IAS/Park
City Mathematics Institute.}
\thanks{
The original manuscript was prepared with the \AmS-\LaTeX\ macro
system and the \Xy-pic\ package.}

\maketitle


\section{Introduction}
\label{sectIntro}
Let $G$ be the group of real points of a semisimple algebraic group defined
over $\QQ$, let $K\subseteq G$ be a maximal compact subgroup, and let
$\G\subset G$ be an arithmetic subgroup.  For every irreducible
representation $(\sigma,U)$ of $G$, Satake \cite{refnSatakeCompact} has
constructed a compactification of the corresponding symmetric space
$D=G/K$.  The procedure is to pass from $(\sigma,U)$ to a representation
$(\pi,V)$ having a nonzero $K$-fixed vector $v$ (a \emph{spherical
representation}); the compactification is then the closure of $D$ under the
embedding $D\hookrightarrow \PP(V)$, $gK \mapsto [\pi(g)v]$.  We denote the
compactification by $\overline{D_\sigma}$ or $\Dstar_\pi$; it is a disjoint
union of \emph{real boundary components}, one of which is $D$ itself and
the others are symmetric spaces of lower rank.  Under certain conditions,
Satake \cite{refnSatakeQuotientCompact} creates a corresponding
compactification $\Xstar_\pi$ of the locally symmetric space $X=\G\back D$.
These conditions were reformulated by Borel
\cite{refnBorelEnsemblesFondamentaux} and christened \emph{geometric
rationality} by Casselman \cite{refnCasselmanGeometricRationality}.

Although the classification of Satake compactifications $\Dstar_\pi$
depends only on $G$ as an $\RR$-group, geometric rationality depends
\emph{a priori} on the $\QQ$-structure.  Borel
\cite[Th\'eor\`eme~4.3]{refnBorelEnsemblesFondamentaux} proved geometric
rationality in the case where $(\sigma,U)$ is strongly rational%
\footnote{This condition is stated in
  \cite[\S1.6]{refnBorelEnsemblesFondamentaux} but only implicitly;
  the term was actually coined later in \cite{refnBorelTits}.
  Unfortunately the adverb ``strongly'' is sometimes omitted in the
  literature, e.g. see \cite{refnCasselmanGeometricRationality}.}
over $\QQ$ in the sense of \cite{refnBorelTits}.  Baily and Borel
\cite{refnBailyBorel} considered the case where $D\subseteq \CC^N$ is a a
bounded symmetric domain and $\Dstar$ is the closure of $D$.  This natural
compactification is topologically equivalent to a certain Satake
compactification \cite[\S5.2]{refnSatakeCompact}, \cite[Theorem
~4]{refnMooreCompactificationsII} and all real boundary components are
Hermitian symmetric spaces.  Baily and Borel prove geometric rationality,
without any rationality condition on $(\sigma,U)$, by a careful
consideration of root systems; another proof was given in
\cite{refnAshMumfordRapoportTai}.  Zucker
\cite{refnZuckerSatakeCompactifications} was the first to raise the general
question of which $(\sigma,U)$ lead to geometrically rational
compactifications.  In \cite{refnCasselmanGeometricRationality} Casselman
gives a criterion for geometric rationality in terms of $\pi$ and the
$\QQ$-index of $G$ \cite{refnTits}.

In this paper we establish geometric rationality in two broad situations.  The
first main result, Theorem ~\ref{thmOriginalBorel} in
\S\ref{sectQrationalRepresentations}, is that $\overline{D_\sigma}$ is
geometrically rational if $(\sigma,U)$ is rational over $\QQ$.  This
generalizes Borel's result \cite{refnBorelEnsemblesFondamentaux} above
which required that $(\sigma,U)$ be strongly rational over $\QQ$.  The
second main result concerns \emph{equal-rank symmetric spaces}, those that
can be written as $G/K$ where $\CCrank G = \rank K$.  Define a \emph{real
equal-rank Satake compactification} to be one for which each real boundary
component is an equal-rank symmetric space; the possibilities have been
listed by
Zucker \cite[(A.2)]{refnZuckerLtwoIHTwo}.  We prove in Theorem
~\ref{thmMainTheorem} of \S\ref{sectEqualRankSatake} that any real
equal-rank Satake compactification is geometrically rational aside from
certain $\QQ$-rank $1$ and $2$ exceptions.  No rationality assumption on
$(\sigma,U)$ or $(\pi,V)$ is made.  Every Hermitian symmetric space is an
equal-rank symmetric space, so this generalizes the geometric rationality
result of Baily and Borel \cite{refnBailyBorel}.  The exceptional cases are
described and treated in Theorem ~\ref{thmSpecialCases} of
\S\ref{sectExceptionalCases}.

The equal-rank result is important since in \cite{refnSaperLModules} the author
proved that the intersection cohomology of a real equal-rank Satake
compactification $\Xstar_\pi$ agrees with both the intersection cohomology
and the weighted cohomology of the reductive Borel-Serre compactification
$\Xhat$.  In the Hermitian case this result on intersection cohomology was
conjectured independently by Rapoport \cite{refnRapoportLetterBorel},
\cite{refnrapoport} and by Goresky and MacPherson
\cite{refnGoreskyMacPhersonWeighted}, while the result on weighted
cohomology in the Hermitian case was proved previously by Goresky, Harder,
and MacPherson \cite{refnGoreskyHarderMacPherson}.

We have chosen to avoid classification theory as much as possible in order
to emphasize the role played by the equal-rank condition; in particular we
do not rely on Zucker's list.  In fact classification theory is only used
in Proposition ~\ref{propFTotallyOrdered} (which uses the classification of
semisimple Lie algebras over $\CC$), in Corollary ~\ref{corFFour} (which
uses the classification of real forms of $F_4$), and in the treatment of
the exceptional cases in Theorem ~\ref{thmSpecialCases}.

Our main results answer several questions raised by Armand Borel
\cite{refnBorelLetterToSaper} in a letter to the author.  It is a pleasure
to thank him and Bill Casselman for helpful correspondence and discussions
concerning this work.  Much of this work took place at the 2002 IAS/Park
City Mathematics Institute on Automorphic Forms; I would like to thank the
organizers of the Institute for providing a stimulating environment.  I
would also like to thank the referee for helpful suggestions regarding the
exposition.

\section{Basic Notions}
\label{sectNotation}
In this paper $G$ will be the group of real points of a semisimple
algebraic $\QQ$-group, $K$ will be a maximal compact subgroup of $G$, and
$D=G/K$ the associated symmetric space.  For a subfield $k$ of $\CC$ we
will find it convenient to speak of a parabolic $k$-subgroup of $G$ when we
strictly mean the group of real points of a parabolic $k$-subgroup of the
algebraic group underlying $G$; similar liberties will be taken for other
subgroups, such as tori or unipotent radicals, and concepts such as roots.
By a representation $(\pi,V)$ of $G$ we will mean a finite-dimensional
$\CC$-vector space $V$ and a homomorphism $\pi\colon G \to \GL(V)$ induced
by a $\CC$-morphism of the underlying algebraic varieties.

Fix a maximal $\QQ$-split torus $\lsb{\QQ}S$ contained in a maximal
$\RR$-split torus $\lsb{\RR} S$ which itself is contained in a maximal
torus $\lsb{\CC}S$ of $G$ defined over $\QQ$.%
\footnote{Much of this paper would become simpler if we could assume that
  $\lsb{\RR}S$ is defined over $\QQ$.  As Borel pointed out to me, this is
  not always possible; there are counterexamples due to Serre.}
Assume compatible orderings on the corresponding root systems have been
chosen.  For $k=\QQ$, $\RR$, or $\CC$, let $\lsb k\D$ denote the simple
$k$-roots and let $\lsb k W$ be the Weyl group of the $k$-root system.  For
$\al\in\lsb k\D$, let $s_\al\in \lsb k W$ denote the corresponding simple
reflection.  The \emph{Coxeter graph} of the system of $k$-roots is a
labeled graph with $\lsb k\D$ for its vertex set and an edge labeled $r$
between $\al\neq \b$ if the product of simple reflections $s_\al s_\beta$
in the Weyl group has order $r>2$.  If we
add an arrow on an edge labeled greater than $3$ which points to the
shorter root we obtain the \emph{Dynkin diagram}.  As customary, we
transfer topological terminology from the graph to $\lsb k\D$.  Thus we may
speak of a connected subset of $\lsb k\D$ or a path in $\lsb k\D$.

The Galois group $\Gal(\CC/k)$ acts on $\lsb\CC\D$ (and hence on
$\CC$-weights) via the $*$-action \cite{refnBorelTits}.  Namely,
$g^*\colon\lsb\CC\D \to \lsb\CC\D$ for $g\in \Gal(\CC/k)$ is defined to be
$w_g\circ g$, where $w_g\in\lsb\CC W$ is uniquely determined by
$w_g(g\lsb\CC\D)=\lsb\CC\D$.  In particular, if $k\subseteq \RR$ we have
$c^*$ where $c$ denotes complex conjugation.  The \emph{opposition
involution} $\iota\colon\lsb\CC\D \to \lsb\CC\D$ is defined similarly by
replacing $g$ with negation.  The $*$-action of $\Gal(\CC/k)$ and $\iota$
commute and are automorphisms of the Dynkin diagram.

Restriction of roots defines
\begin{equation*}
\restr_{\CC/k}\colon\lsb\CC\D \to \lsb k\D \cup \{0\}.
\end{equation*}
The fibers $\restr_{\CC/k}^{-1}(\al)$ for $\al \in \lsb k\D$ are nonempty
Galois orbits; the fiber $\D^0_{\CC/k} = \restr_{\CC/k}^{-1}(0)$ is Galois
invariant and its elements are the \emph{$k$-anisotropic roots}.  For
$\theta\subseteq \lsb k\D$ define $\e_{\CC/k}(\theta) =
\restr_{\CC/k}^{-1}(\theta \cup\{0\})$.  These are the types of the
parabolic $k$-subgroups of $G$.

The \emph{$k$-index} of $G$ \cite{refnTits} consists of the Dynkin diagram
for $\lsb\CC\D$, the subset $\D^0_{\CC/k}$, and the $*$-action of
$\Gal(\CC/k)$.  This can be represented diagrammatically; in the case
$k=\RR$ one recovers the \emph{Satake diagram}.  We will often by abuse of
notation refer to the index simply as $\lsb\CC\D$.

We will use the following repeatedly \cite[Proposition~6.15]{refnBorelTits}:
\begin{prop}
\label{propConnected}
For $\al$, $\b\in\lsb k\D$, $\{\al,\b\}$ is connected if and only if for
every $\tilde\al\in \restr_{\CC/k}^{-1}(\al)$, there exists $\tilde\b\in
\restr_{\CC/k}^{-1}(\b)$ and $\psi\subseteq \D^0_{\CC/k}$ such that
$\psi\cup\{\tilde\al,\tilde\b\}$ is connected.
\end{prop}


Let $\chi$ be the highest weight of an irreducible representation
$(\pi,V)$ of $G$ and set $\d = \{\,\al\in\lsb \CC\D\mid s_\al\chi \neq
\chi\,\}$.  (If necessary we will denote this $\d_\pi$ or $\d_\chi$.)  Let
$\lsb k\chi = \chi|_{\lsb k S}$ and define $\lsb k\d\subseteq \lsb k \D$
analogously.  An analogue of Proposition ~\ref{propConnected} shows that
$\b\in\lsb k\d$ if and only if there exists $\tilde\b\in
\restr_{\CC/k}^{-1}(\b)$ and $\psi\subseteq \D^0_{\CC/k}$ such that
$\psi\cup\{\tilde\b\}$ is connected and contains an element of $\d$
\cite[Corollary~7.2]{refnCasselmanGeometricRationality}, \cite[Remark in
(2.4)]{refnZuckerSatakeCompactifications}.

For a linear combination of simple $k$-roots $\sum d_\alpha \alpha$, the
\emph{support} $\supp_k(\sum d_\alpha \alpha)$ is defined as
$\{\alpha\in\lsb k\D\mid d_\alpha\neq 0\}$; the linear combination is
called \emph{codominant} if all $d_\alpha\ge0 $.  A subset $\theta\subseteq
\lsb k\D$ is called \emph{$\d$\dash connected} if every connected component
of $\theta$ contains an element of $\lsb k\d$.  For every $k$-weight $\lambda$
of $(\pi,V)$, the difference $\lsb k\chi - \lambda$ is codominant with
$\d$-connected support (and integral coefficients).  In fact every
$\d$-connected subset of $\lsb k \D$ arises in this way
\cite[12.16]{refnBorelTits}.

For $\theta\subseteq \lsb k\D$, let $\kap(\theta)$ denote the largest
$\d$-connected subset of $\theta$ and let $\o(\theta)$ denote the largest
subset $\Upsilon$ of $\lsb k\D$ with $\kap(\Upsilon)=\kap(\theta)$.
Clearly
\begin{equation*}
\kap(\theta)\subseteq \theta \subseteq \o(\theta).
\end{equation*}
Let $\z(\theta)$ denote the complement of $\kap(\theta)$ in $\o(\theta)$.
Equivalently, $\z(\theta)$ consists of those roots that are not in $\lsb
k\d$, not in $\kap(\theta)$ and not joined by an edge to a root in
$\kap(\theta)$.  It will sometimes be useful, especially in
\S\ref{sectEqualRankSatake}, to use the notation
\begin{equation*}
\theta^+ = \theta \cup \{\,\al\in \lsb k\D\mid \text{$\al$ is connected by
  an edge to a root in $\theta$}\,\}.
\end{equation*}

We need to recall the various notions of rationality for a representation
\cite[\S12]{refnBorelTits}.  Assume that $V$ has a $k$-structure.  An irreducible representation
$(\pi,V)$ is called \emph{projectively rational over $k$} if the associated
projective representation $\pi'\colon G\to \PGL(V)$ is defined over $k$.
It is called \emph{rational over $k$} if $\pi\colon G \to \GL(V)$ is itself
defined over $k$; it is called \emph{strongly rational over $k$} if
furthermore the
parabolic subgroup $P_\pi$ stabilizing the line spanned by the highest
weight vector is defined over $k$.  In terms of the highest weight $\chi$
of $(\pi,V)$, the representation is projectively rational over $k$ if and
only if $\chi$ is $\Gal(\CC/k)$-invariant under the $*$-action and it is
strongly rational over $k$ if and only if
in addition $s_\al\chi=\chi$ for all $\al\in\D^0_{\CC/k}$.

An irreducible representation $(\pi,V)$ is called \emph{spherical}
if there exists a nonzero $K$-fixed vector $v\in V$; such a vector $v$ is
unique up to scalar multiplication.  By a theorem of Helgason
\cite{refnHelgasonDuality}, \cite[Theorem~8.49]{refnKnapp}, $(\pi,V)$ is
spherical if and only if its highest weight $\chi$ is trivial on the
maximal $\RR$-anisotropic subtorus of $\lsb{\CC}S$ and
$(\chi,\alpha)/(\alpha,\alpha)\in \ZZ$ for all $\alpha \in \lsb{\RR}\D$.
In particular, a spherical representation $(\pi,V)$ is strongly rational
over $\RR$.  In this case $\d$ is $\Gal(\CC/\RR)$-invariant and
disjoint from $\D_{\CC/\RR}^0$; it follows that
$\epsilon_{\CC/\RR}(\o(\theta)) = \o(\epsilon_{\CC/\RR}(\theta))$ for all
$\theta \subseteq \lsb\RR\D$.  It also follows from \cite[proof of
Theorem~8.49]{refnKnapp} that a $K$-fixed vector $v$ has a nonzero
component along any highest weight vector.

\section{Satake Compactifications}
\label{sectSatakeCompact}
As in \cite{refnCasselmanGeometricRationality} we take as our starting
point an irreducible spherical representation $(\pi,V)$ of $G$,
nontrivial on each noncompact $\RR$-simple factor of $G$, with $K$-fixed
vector $v$; for the relation with Satake's original construction, see
\S\ref{sectSatakeCompactOriginal}.  The Satake compactification
$\Dstar_\pi$ associated to $(\pi,V)$ is defined to be the closure of the
image of $D=G/K$ under the embedding $D\hookrightarrow \PP(V)$, $gK \mapsto
[\pi(g)v]$.  The action of $G$ on $D$ extends to an action on $\Dstar_\pi$.
For every parabolic $\RR$-subgroup $P$, the subset $D_{P,h} \subseteq
\Dstar_\pi$ of points fixed by $N_P$, the unipotent radical, is called a
\emph{real boundary component}.  The Satake compactification is the
disjoint union of the real boundary components, however different $P$ may
yield the same real boundary component.  By associating to each real
boundary component its normalizer we obtain a one-to-one correspondence.

A parabolic $\RR$-subgroup $P$ is \emph{$\d$-saturated} if it is conjugate
to a standard parabolic $\RR$-subgroup with type
$\epsilon_{\CC/\RR}(\o(\theta))$ for some $\theta\subseteq \lsb\RR\D$.  The
subgroups that arise as normalizers of real boundary components are the
precisely the $\d$-saturated parabolic $\RR$-subgroups.  The action of $P$
on $D_{P,h}$ descends to an action of its Levi quotient $L_P=P/N_P$ and the
subgroup $L_{P,\l}\subseteq L_P$ which fixes $D_{P,h}$ pointwise will be
called the \emph{centralizer group}.  It has as its identity
component (in the Zariski topology) the maximal normal connected
$\RR$-subgroup with simple $\RR$-roots $\z(\theta)$; the simple $\CC$-roots
of $L_{P,\l}$ are $\z(\e_{\CC/\RR}(\theta))$.  Thus $D_{P,h}$ is the
symmetric space corresponding to $L_{P,h}= L_P/L_{P,\l}$, a semisimple
$\RR$-group with simple $\RR$-roots $\kap(\theta)$; the simple $\CC$-roots
of $L_{P,h}$ are $\kap(\e_{\CC/\RR}(\theta))$, the connected components of
$\e_{\CC/\RR}(\kap(\theta))$ which are not wholly contained in
$\D^0_{\CC/\RR}$.

If $D_{P,h}$ and $D_{P',h}$ are two standard real boundary components with
normalizers having type $\epsilon_{\CC/\RR}(\o(\theta))$ and
$\epsilon_{\CC/\RR}(\o(\theta'))$ respectively, define $D_{P,h} \le
D_{P',h}$ if $\kap(\theta)\subseteq \kap(\theta')$; this is a partial order
on the standard real boundary components.

More generally the same procedure allows one to associate a Satake
compactification $\Dstar_\pi$ to a triple $(\pi,V,v)$ consisting of a (not
necessarily irreducible) representation $(\pi,V)$ which is
nontrivial on each noncompact $\RR$-simple factor of $G$ and a $K$-fixed
vector $v\in V$ whose $G$-orbit spans $V$.  Note that every irreducible
component of $(\pi,V)$ is automatically spherical.  We will only consider
$(\pi,V,v)$ satisfying a very restrictive condition (although Casselman
points out in \cite{refnCasselmanGeometricRationality} that the general
case is worthy of further study):

(R) Assume that $(\pi,V)$ may be decomposed into the direct sum of
submodules $(\pi_0,V_0)$ and $(\pi', V')$ such that $(\pi_0,V_0)$ is
irreducible with highest weight $\chi_0$ and the difference $\chi_0-\chi'$,
where $\chi'$ is the highest weight of any irreducible component of
$(\pi',V')$, is codominant with $\d_{\chi_0}$-connected support.

\begin{lem}
\label{lemReducibleSatakeCompact}
If $(\pi,V,v)$ satisfies condition \textup{(R)}, then
$\Dstar_\pi\cong \Dstar_{\pi_0}$ as $G$-spaces.
\end{lem}
\begin{proof}
Let $y\in \Dstar_\pi \setminus D$ and let $y_i = [\pi(g_i)v]\in D$ be a
sequence converging to $y$.  Write $g_iK = k_i a_i K \in K\overline{A^+}K$
by the Cartan decomposition \cite[IX, Theorem 1.1]{refnHelgason}; here $A$
is the identity component (in the classical topology) of $\lsb{\RR}S$ and
$\overline{A^+}= \{\,a\in A\mid a^\alpha\ge 1,\ \alpha\in \lsb{\RR}\D\,\}$.
By passing to a subsequence and conjugating $\lsb{\RR}S$ if necessary, we
may assume that $k_i \to e$ in which case we may assume simply that $y_i =
[\pi(a_i)v]$.  Choose homogeneous coordinates $x_\lambda$ on $\PP(V)$
corresponding to a basis of $\CC$-weight vectors of $V$.  We know there
exists a $\CC$-weight $\lambda$ of $V$ such that $y$ is not in the zero set
of
$x_\lambda$; let $\chi'$ be the highest weight of an irreducible
component of $V$ containing the $\CC$-weight $\lambda$.  Then
$({x_{\chi_0}}/{x_\lambda})(y_i) \propto a_i^{\chi_0-\chi'+\chi'- \lambda}
\ge 1$ since the exponent is codominant; consequently $y$ is not in the
zero set of $x_{\chi_0}$.

Since  $\Dstar_\pi\cap \PP(V')$ is therefore empty, the 
projection $V\to V/V'\cong V_0$ induces a continuous map $\Dstar_\pi \to
\Dstar_{\pi_0}\subseteq \PP(V_0)$ of $G$-spaces.
Write $v=v_0 +
v'$.  In order to see this map is bijective, we need to show that if a
sequence $y_{0,i} = [\pi_0(a_i)v_0]$ converges to $y_0\in \Dstar_{\pi_0}$,
then the sequence $y_i = [\pi(a_i)v]$ converges in $\Dstar_\pi$.  By
hypothesis,
$({x_{\lambda_0}}/{x_{\chi_0}})(y_{0,i}) \propto a_i^{-(\chi_0-\lambda_0)}$
converges for every $\CC$-weight $\lambda_0$ of $V_0$.
Consequently for any  $\d_{\chi_0}$-connected subset $\theta_0\subseteq
\lsb{\CC}\D$, either $a_i^{-\al}$ converges for all $\al\in\theta_0$ or
there exists $\al\in\theta_0$ with $a_i^{-\al}\to0$; the proof is by induction
on $\#\theta_0$.  Since $a_i^{-\al}\le1$ for all $\al$, this
implies
\begin{equation}
\label{eqnConvergenceCondition}
\text{$a_i^{-\nu}$ converges for all $\nu$ codominant with
$\d_{\chi_0}$-connected support.}
\end{equation}
We need to show for any $\CC$-weight $\lambda$ of an irreducible component
with highest weight $\chi'$ that $({x_\lambda}/{x_{\chi_0}})(y_i) \propto
a_i^{-(\chi_0-\chi'+\chi'- \lambda)}$ converges.  Since $\supp(\chi'-
\lambda)$ is $\d_{\chi'}$-connected and $\d_{\chi'} \subseteq \d_{\chi_0}
\cup \supp(\chi_0-\chi')$, it follows from (R) that
$\supp(\chi_0-\chi')\cup \supp(\chi'-\lambda)$ is $\d_{\chi_0}$-connected.
Convergence then follows from \eqref{eqnConvergenceCondition}.
\end{proof}

\section{Satake's Original Construction}
\label{sectSatakeCompactOriginal}
In this section we relate Satake's original construction
\cite{refnSatakeCompact} of compactifications of $D$ to those considered in
\S\ref{sectSatakeCompact}.  Let $(\sigma,U)$ be an irreducible
representation of $G$ (not necessarily spherical) which is nontrivial on
each $\RR$-simple factor.  Fix a positive definite Hermitian inner product
on $U$ which is \emph{admissible} in the sense that $ \sigma(g)^* =
\sigma(\theta g)^{-1}$ (here $\theta$ denotes the Cartan involution
associated to $K$); such an inner product always exists
\cite{refnMatsushimaMurakami} and is unique up to rescaling.  Let $S(U)$
denote the corresponding real vector space of self-adjoint endomorphisms of
$U$.  Satake defines the compactification $\overline{D_\sigma}$ associated
to $(\sigma,U)$ to be the closure of the image of $D$ under the embedding
$D\hookrightarrow \PP S(U)$, $gK \mapsto \sigma(g)\sigma(g)^*$.

Let $g\in G$ act on $S(U)$ by $T\mapsto \sigma(g) \circ T \circ
\sigma(g)^*= \sigma(g) \circ T \circ \sigma(\theta g)^{-1}$ and extend to a
representation on $S(U)_{\CC}$, the space of all endomorphisms of
$U$.  This representation is $U \otimes \lsp\theta U^*$, where $\lsp\theta
U^*$ is the twist by $\theta$ of the usual contragredient representation.
The identity endomorphism $I\in U \otimes \lsp\theta U^*$ is a $K$-fixed
vector and we let $(\pi,V)$ denote the smallest subrepresentation
containing $I$.  It is clear that $\overline{D_\sigma}$ is isomorphic
with the Satake compactification $\Dstar_\pi$ associated to
$(\pi,V,I)$ as in \S\ref{sectSatakeCompact}.

\begin{prop}
\label{propConditionR}
The triple $(\pi,V,I)$ satisfies condition \textup{(R)} with $\chi_0 = 2\Re
\mu$, where $\mu$ is the highest weight of $(\sigma, U)$.
\end{prop}
\begin{proof}
The $\CC$-weights of the $K$-fixed vector $I$ include the highest weight
$\chi'$ of any irreducible component of $(\pi,V)$.  However $I=
\sum_\lambda u_\lambda\otimes u_\lambda^*$, where $u_\lambda$ ranges over a
basis of $\CC$-weight vectors of $U$ and $u_\lambda^*$ ranges over the
corresponding dual basis.  Thus $\chi' = \lambda -\theta \lambda = 2\Re
\lambda$, where $\lambda$ is a $\CC$-weight of $U$.  Since $\Re\lambda$ is
an $\RR$-weight of $U$, the difference $\Re\mu - \Re \lambda$ is
$\RR$-codominant
and $\supp_\RR(\Re\mu - \Re \lambda)$ is $\lsb{\RR}\d_\mu$-connected.  If
$\alpha\in\lsb\RR\D$, then $\supp_\CC(2\alpha)$ has at most $2$ connected
components, each containing an element of $\restr_{\CC/\RR}^{-1}(\alpha)$.
We see that $\supp_{\CC}(2\Re\mu -2\Re\lambda)$ is $\d_{2\Re\mu}$-connected
since $\d_{2\Re\mu} = \restr_{\CC/\RR}^{-1}(\lsb{\RR}\d_\mu)$ and since if
$(\alpha,\beta)\neq 0$ then each connected component of
$\supp_\CC(2\alpha)$ is connected to a connected component of
$\supp_{\CC}(2\beta)$ \cite[(1.7)]{refnZuckerSatakeCompactifications}.
Finally a component $(\pi_0,V_0)$ with highest weight $2\Re \mu$ does in
fact occur, since one may check that $\smash[b]{\sum_{\lambda|_{\lsb{\RR}S}
=\mu|_{\lsb{\RR}S}} u_\lambda\otimes u_\lambda^*}$ is a highest weight
vector.%
\footnote{A more direct proof that $2\Re \mu$ is the highest weight of an
  irreducible spherical representation is given by Harish-Chandra
  \cite[Lemma 2]{refnHarishChandraSphericalI}.}
\end{proof}
\begin{cor}
\label{corIdentifyOriginalSatake}
If $(\sigma,U)$ has highest weight $\mu$ and $(\pi_0,V_0)$ is the
irreducible spherical representation with highest weight $2\Re \mu$, then
$\overline{D_\mu}\cong \Dstar_{\pi_0}$ as $G$-spaces.
\end{cor}

\section{Geometric Rationality}
A real boundary component $D_{P,h}$ of $\Dstar_\pi$ is called
{\itshape rational\/} \cite[\S\S3.5,~3.6]{refnBailyBorel} if
\begin{enumerate} 
\item its normalizer $P$ is defined over $\QQ$, and
\label{itemRationalNormalizer}
\item the centralizer group $L_{P,\l}$ contains a normal subgroup $\tilde
L_{P,\l}$ of $L_P$ defined over $\QQ$ such that $L_{P,\l}/\tilde L_{R,\l}$
is compact.
\label{itemRationalCentralizer}
\end{enumerate}
A Satake compactification $\Dstar_\pi$ is called
\emph{geometrically rational} if every real boundary component $D_{P,h}$
whose normalizer has type $\epsilon_{\CC/\RR}(\o(\epsilon_{\RR/\QQ}(\U)))$
for some $\d$-connected subset $\U\subseteq \lsb\QQ\D$ is rational.  If
$\QQrank G=0$ (and hence $X=\G\back D$ is already compact) any Satake
compactification $\Dstar_\pi$ is geometrically rational, so our
main interest is when $\QQrank G>0$.

Casselman \cite[Theorems ~8.2, 8.4]{refnCasselmanGeometricRationality}
proves the following criterion:

\begin{thm}
\label{thmCasselmanCriterion}
The Satake compactification associated to $(\pi,V)$ is geometrically
rational if and only if
\begin{enumerate}
\item $\o(\D_{\CC/\QQ}^0)$  is Galois invariant, and
\label{itemGRone}
\item $\kap(\D_{\CC/\QQ}^0)$ is Galois invariant modulo $\RR$-anisotropic
  roots, that is, the \linebreak $\Gal(\CC/\QQ)$-orbit of
  $\kap(\D_{\CC/\QQ}^0)$ is contained in $\kap(\D_{\CC/\QQ}^0) \cup
  \D_{\CC/\RR}^0$.
\label{itemGRtwo}
\end{enumerate}
\end{thm}

\section{Representations Rational over $\QQ$}
\label{sectQrationalRepresentations}
\begin{thm}
\label{thmDeltaRational}
The Satake compactification associated to $(\pi,V)$ is geometrically
rational if $\d$ is Galois invariant.
\end{thm}

\begin{proof}
We need to verify \itemref{itemGRone} and \itemref{itemGRtwo} from Theorem
~\ref{thmCasselmanCriterion}.  Now $\kap(\D_{\CC/\QQ}^0)$ is the union of
those connected components of $\D_{\CC/\QQ}^0$ which contain an element of
$\d$.  Since both $\d$ and $\D_{\CC/\QQ}^0$ are Galois invariant, clearly
$\kap(\D_{\CC/\QQ}^0)$ is Galois invariant.  This proves
\itemref{itemGRtwo} (in fact $\RR$-anisotropic roots are not needed).  Now
$\z(\D_{\CC/\QQ}^0)$ consists of roots not in $\d$ and not in
$\kap(\D^0_{\CC/\QQ})^+$, hence it is Galois invariant.  Since
$\o(\D_{\CC/\QQ}^0) = \kap(\D_{\CC/\QQ}^0) \cup \z(\D_{\CC/\QQ}^0)$ this
proves \itemref{itemGRone}.
\end{proof}

\begin{cor}
\label{corBorel}
The Satake compactification associated to $(\pi,V)$ is geometrically
rational if $(\pi,V)$ is projectively rational over $\QQ$.
\end{cor}

\begin{proof}
By \cite[\S12.6]{refnBorelTits}, the associated projective representation
$\pi'\colon G \to \PGL(V)$ is defined over $\QQ$ if and only if the highest
weight of $V$ is Galois invariant; thus $\d$ is Galois invariant.
\end{proof}

Note that in Satake's original construction, the passage from $(\sigma,U)$
to $(\pi_0,V_0)$ does not necessarily preserve rationality.  Thus Corollary
~\ref{corBorel} does not settle the question of geometric rationality if
$(\sigma,U)$ is assumed to be rational over $\QQ$.  Nonetheless we can
prove:%
\footnote{The result of Theorem ~\ref{thmOriginalBorel} is asserted in
  \cite[(3.3), (3.4)]{refnZuckerSatakeCompactifications}, however in
  \cite[\S9]{refnCasselmanGeometricRationality} it is noted that
  \cite[Prop.~(3.3)(ii)]{refnZuckerSatakeCompactifications} is incorrect.
  Also the proof of \cite[Prop.~(3.4)]{refnZuckerSatakeCompactifications}
  seems to implicitly assume that $(\sigma,U)$ is strongly rational over
  $\RR$.}
\begin{thm}
\label{thmOriginalBorel}
The Satake compactification associated to $(\sigma,U)$ by Satake's original
construction is geometrically rational if $(\sigma,U)$ is projectively
rational over $\QQ$.
\end{thm}
\begin{proof}
By Corollary ~\ref{corIdentifyOriginalSatake}, we need to verify
\itemref{itemGRone} and \itemref{itemGRtwo} from Theorem
~\ref{thmCasselmanCriterion} for $\d = \d_{2\Re \mu}$ under the hypothesis
that $\d_\mu$ is Galois invariant.  By the remarks following Proposition
~\ref{propConnected} and the fact that $\d_\mu$ is $c^*$-invariant, we can
write
\begin{equation*}
\d= \{\,\alpha\in
\lsb\CC\D\setminus \D_{\CC/\RR}^0 \mid \text{there exists $\psi\subseteq
  \D_{\CC/\RR}^0$ such that $\psi\cup\{\alpha\}$ is $\d_\mu$-connected}\,\}.
\end{equation*}
Let $\kap=\kap(\D_{\CC/\QQ}^0)$ is the union of those connected components
of $\D_{\CC/\QQ}^0$ which contain an element of $\d$; similarly let
$\kap_\mu= \kap_\mu(\D_{\CC/\QQ}^0)$ is the union of those connected
components of $\D_{\CC/\QQ}^0$ which contain an element of $\d_\mu$.
Clearly $\kap\subseteq \kap_\mu$ while $\kap_\mu \setminus \kap$ is the
union of those connected components of $\D_{\CC/\QQ}^0$ which lie within
$\D_{\CC/\RR}^0$ and contain an element of $\d_\mu$.  Since $\kap_\mu$ is
Galois invariant, \itemref{itemGRtwo} is satisfied.

For \itemref{itemGRone} we calculate
$\lsb\CC\D\setminus \o(\D_{\CC/\QQ}^0)=(\lsb\CC\D\setminus\kap) \cap (\d
\cup \kap^+) = (\lsb\CC\D\setminus\kap_\mu) \cap (\d \cup \kap_\mu^+)$,
where the second equality holds since $\kap_\mu\setminus \kap\subseteq
\D_{\CC/\RR}^0$ is disjoint from $\d \cup \kap^+$ and since
$(\kap_\mu^+\setminus \kap^+)\setminus \kap_\mu \subseteq \d$.  Since
$\d_\mu\cap (\lsb\CC\D\setminus \D_{\CC/\QQ}^0) \subseteq \d \subseteq
(\d_\mu\cap (\lsb\CC\D\setminus \D_{\CC/\QQ}^0)) \cup \kap_\mu^+$, we can
thus write
$\lsb\CC\D\setminus \o(\D_{\CC/\QQ}^0)= (\lsb\CC\D\setminus\kap_\mu) \cap
((\d_\mu\cap (\lsb\CC\D\setminus \D_{\CC/\QQ}^0)) \cup \kap_\mu^+)$,
which is clearly Galois invariant.
\end{proof}

\section{Equal-rank Satake Compactifications}
\label{sectEqualRankSatake}
We begin with some generalities about the structure of root systems, here
applied to the simple $\CC$-roots $\lsb\CC\D$ of $G$.  Let $\iota$ be the
opposition involution, that is, the negative of the longest element of the
Weyl group.  For a connected subset $\psi\subseteq \lsb\CC\D$
which is invariant under $\iota$, let $\iota|_\psi$ denote the restriction,
while $\iota_\psi$ denotes the opposition involution of the subroot system
with simple roots $\psi$.  Recall
\begin{equation*}
\psi^+ = \psi \cup \{\,\al\in \lsb\CC\D\mid \text{$\al$ is connected by an
  edge to a root in $\psi$}\,\}
\end{equation*}

\begin{defn}
\label{defnF}
Let $\Ftilde$ denote the family of nonempty connected $\iota$-invariant
subsets $\psi\subseteq \lsb\CC\D$ for which $\iota|_\psi=\iota_\psi$.  Let
$\F$ consist of those $\psi\in \Ftilde$ such that
\begin{enumerate}
\item \label{itemFone} $\psi^+\setminus \psi$ modulo $\iota$ has
  cardinality $\le 1$, and
\item \label{itemFtwo} there exists $\psi'\supseteq \psi^+$,
  $\psi'\in\Ftilde$ satisfying \itemref{itemFone}, such that all components
  of $\psi'\setminus \psi^+$ are in $\Ftilde$.
\end{enumerate}
Let $\Ftilde^*\subseteq \Ftilde$ exclude sets of cardinality $1$ which are
not components of $\lsb\CC\D$ and for each component $C\subseteq
\lsb\CC\D$, let $\Ftilde_C = \{\psi\in\Ftilde\mid \psi\subseteq C\}$;
similarly define $\F^*$ and $\F_C$.
\end{defn}

We view $\F$ and $\Ftilde$ as partially ordered sets by inclusion.  For any
partially ordered set $\mathscr{P}$, recall that $\psi'\in\mathscr{P}$
\emph{covers} $\psi\in\mathscr{P}$ if $\psi < \psi'$ and there does not
exists $\psi''\in\mathscr{P}$ with $\psi<\psi''< \psi'$.  The \emph{Hasse
diagram} of the partially ordered set $\mathscr{P}$ has nodes $\mathscr{P}$
and edges the cover relations; if $\psi'$ covers $\psi$ we draw $\psi'$ to
the right of $\psi$.  (For the basic terminology of partially ordered sets,
see \cite[Chapter~3]{refnStanleyEnumCombI}, but note that in
\cite{refnStanleyEnumCombI} $\psi'$ is drawn \emph{above} $\psi$.)

\begin{prop}
\label{propFTotallyOrdered}
The Hasse diagrams of $\F_C$ and $\Ftilde_C$ are given in Figure
~\textup{\ref{figHasseDiagrams}}.  In particular if $\psi<\psi'\in\Ftilde$
and $\psi\in\F$ then $\psi'\in\F$.  Also if $C$ is not type $F_4$,
then\textup:
\begin{enumerate}
\item
\label{itemNotFFour}
$\Ftilde^*_C$ is totally ordered.
\item 
\label{itemIncomparable}
Suppose $\psi$, $\psi'\in\F_C$ are incomparable.
If $C$ is not type $B_n$, $C_n$, or $G_2$, then one of $\psi$ and $\psi'$
is type $A_1$ and $\psi\cup\psi'$ is disconnected.  If $C$ is type $B_n$,
$C_n$, or $G_2$, one of $\psi$ and $\psi'$ is type $A_1$ and is covered by
$C$.

\end{enumerate}
\end{prop}
\begin{proof}
Recall \cite{refnTits} that $\iota|_C$ is the unique nontrivial involution
in the cases $A_n$ ($n>1$), $D_n$ ($n>4$, odd), and $E_6$, and is trivial
otherwise.  If $\rank C>1$, let $C_0\subseteq C$ be defined as follows: if
$\iota|_C$ is nontrivial, $C_0$ is the unique subdiagram of type $A_2$ or
$A_3$ which is nontrivially stabilized under $\iota|_C$; if $\iota|_C$ is
trivial, $C_0$ is the unique subdiagram of type $B_2$, $D_4$, or $G_2$.  The
subdiagram $C_0$ ``determines'' $\iota|_C$.  In particular if
$\psi\subseteq C$ is a
connected $\iota$-invariant subset with $\rank \psi > 1$, then
$\iota|_\psi=\iota_\psi$ is equivalent to $C_0 = \psi_0$.  Such $\psi$ are
easy to enumerate.  For example, if $C$ is type $B_n$, then $\psi$ can be
any connected segment containing the double bond at one end.
Together with
the cardinality $1$ subsets, such calculations yield $\Ftilde_C$.  To
determine $\F_C$, one then checks which of these $\psi$ satisfy the
additional conditions of Definition ~\ref{defnF}.  The results are pictured
in Figure ~\ref{figHasseDiagrams} and the proposition follows.
\begin{figure}
\begin{align*}
&\vcenter{\xymatrix @!0 @M=0pt @R-.1in {
{\rlap{\text{(Type $A_n$, $n$ odd)}}}&{}&{}&{}&{}&{\bullet}\save[]+<0in,-.1in>*{\scriptstyle A_1}\restore \ar@{-}'[r] '[rr] [rrr] &
{\bullet}\save[]+<0in,-.1in>*{\scriptstyle A_3}\restore &
{\bullet}\save[]+<0in,-.1in>*{\scriptstyle A_5}\restore &
{} \ar@{}[r]|{.\ .\ .} &
{} \ar@{-}'[r] [rr] &
{\bullet}\save[]+<0in,-.1in>*{\scriptstyle A_{n-2}}\restore &
{\bullet}\save[]+<0in,-.1in>*{\scriptstyle A_n}\restore
}} \\
&\vcenter{\xymatrix @!0 @M=0pt @R-.1in {
{\rlap{\text{(Type $A_n$, $n$ even)}}}&{}&{}&{}&{}&{\bullet}\save[]+<0in,-.1in>*{\scriptstyle A_2}\restore \ar@{-}'[r]  '[rr] [rrr] &
{\bullet}\save[]+<0in,-.1in>*{\scriptstyle A_4}\restore &
{\bullet}\save[]+<0in,-.1in>*{\scriptstyle A_6}\restore &
{} \ar@{}[r]|{.\ .\ .} &
{} \ar@{-}'[r] [rr] &
{\bullet}\save[]+<0in,-.1in>*{\scriptstyle A_{n-2}}\restore &
{\bullet}\save[]+<0in,-.1in>*{\scriptstyle A_n}\restore
}} \\
&\vcenter{\xymatrix @!0 @M=0pt @R-.1in {
&{}{}&{}&{}&{}&{}&{}&{}&{}&{}&{\bullet}\save[]+<0in,-.1in>*{\scriptstyle A_{1}}\restore \ar@{-}[dr]\\
{\rlap{\text{(Type $B_n$)}}}&{}&{}&{}&{}&{\bullet}\save[]+<0in,-.1in>*{\scriptstyle A_1}\restore \ar@{-}'[r] '[rr] [rrr] &
{\bullet}\save[]+<0in,-.1in>*{\scriptstyle B_2}\restore &
{\bullet}\save[]+<0in,-.1in>*{\scriptstyle B_3}\restore &
{} \ar@{}[r]|{.\ .\ .} &
{} \ar@{-}'[r] [rr] &
{\bullet}\save[]+<0in,-.1in>*{\scriptstyle B_{n-1}}\restore &
{\bullet}\save[]+<0in,-.1in>*{\scriptstyle B_n}\restore
}} \\
&\vcenter{\xymatrix @!0 @M=0pt @R-.1in {
&{}{}&{}&{}&{}&{}&{}&{}&{}&{}&{\bullet}\save[]+<0in,-.1in>*{\scriptstyle A_{1}}\restore \ar@{-}[dr]\\
{\rlap{\text{(Type $C_n$)}}}&{}&{}&{}&{}&{\bullet}\save[]+<0in,-.1in>*{\scriptstyle A_1}\restore \ar@{-}'[r] '[rr] [rrr] &
{\bullet}\save[]+<0in,-.1in>*{\scriptstyle C_2}\restore &
{\bullet}\save[]+<0in,-.1in>*{\scriptstyle C_3}\restore &
{} \ar@{}[r]|{.\ .\ .} &
{} \ar@{-}'[r] [rr] &
{\bullet}\save[]+<0in,-.1in>*{\scriptstyle C_{n-1}}\restore &
{\bullet}\save[]+<0in,-.1in>*{\scriptstyle C_n}\restore
}} \\
&\vcenter{\xymatrix @!0 @M=0pt @R-.1in {
&{}{}&{}&{}&{}&{}&{}&{}&{}&{}&{\bullet}\save[]+<0in,-.1in>*{\scriptstyle
A_{1}}\restore \ar@{-}[dr]\\
{\rlap{\text{(Type $D_n$, $n$ odd)}}}&{}&{}&{}&{}&{\bullet}\save[]+<0in,-.1in>*{\scriptstyle A_3}\restore \ar@{-}'[r] '[rr] [rrr] &
{\bullet}\save[]+<0in,-.1in>*{\scriptstyle D_5}\restore &
{\bullet}\save[]+<0in,-.1in>*{\scriptstyle D_7}\restore &
{} \ar@{}[r]|{.\ .\ .} &
{} \ar@{-}'[r] [rr] &
{\bullet}\save[]+<0in,-.1in>*{\scriptstyle D_{n-2}}\restore &
{\bullet}\save[]+<0in,-.1in>*{\scriptstyle D_n}\restore \\
}} \\
&\vcenter{\xymatrix @!0 @M=0pt @R-.1in{
&{}&{}&{}&{}&{\bullet}\save[]+<0in,-.1in>*{\scriptstyle A_{1}}\restore \ar@{-}[dr] &
{}&{}&{}&{}&{\bullet}\save[]+<0in,-.1in>*{\scriptstyle A_{1}}\restore \ar@{-}[dr]\\
{\rlap{\text{(Type $D_n$, $n>4$ even)}}}&{}&{}&{}&{}&{}&{\bullet}\save[]+<0in,-.1in>*{\scriptstyle D_4}\restore  \ar@{-}'[r] [rr] &
{\bullet}\save[]+<0in,-.1in>*{\scriptstyle D_6}\restore &
{} \ar@{}[r]|{.\ .\ .} &
{} \ar@{-}'[r] [rr] &
{\bullet}\save[]+<0in,-.1in>*{\scriptstyle D_{n-2}}\restore &
{\bullet}\save[]+<0in,-.1in>*{\scriptstyle D_n}\restore \\
&{}&{}&{}&{}&{\bullet}\save[]+<0in,-.1in>*{\scriptstyle A_{1}}\restore \ar@{-}[ur]
}} \\
&\vcenter{\xymatrix @!0 @M=0pt @R-.1in{
&{}&{}&{}&{}&
{\bullet}\save[]+<0in,-.1in>*{\scriptstyle A_{1}}\restore \ar@{-}[dr]\\
{\rlap{\text{(Type $D_4$)}}}&{}&{}&{}&{}&
{\bullet}\save[]+<0in,-.1in>*{\scriptstyle A_{1}}\restore \ar@{-}[r]&
{\bullet}\save[]+<0in,-.1in>*{\scriptstyle D_4}\restore \\
&{}&{}&{}&{}&{\bullet}\save[]+<0in,-.1in>*{\scriptstyle A_{1}}\restore \ar@{-}[ur]
}} \\
&\vcenter{\xymatrix @!0 @M=0pt @R-.1in {
{\rlap{\text{(Type $E_6$)}}}&{}&{}&{}&{}&{\circ}\save[]+<0in,-.1in>*{\scriptstyle A_3}\restore \ar@{.}[r] &
{\bullet}\save[]+<0in,-.1in>*{\scriptstyle A_5}\restore  \ar@{-}[r]&
{\bullet}\save[]+<0in,-.1in>*{\scriptstyle E_6}\restore
}} \\
&\vcenter{\xymatrix @!0 @M=0pt @R-.1in {
&{}&{}&{}&{}&{\circ}\save[]+<0in,-.1in>*{\scriptstyle D_4}\restore
      \ar@{.}[dr]  \\
{\rlap{\text{(Type $E_7$)}}}&{}&{}&{}&{}&
{\bullet}\save[]+<0in,-.1in>*{\scriptstyle A_1}\restore \ar@{-}'[r] [rr]&
{\bullet}\save[]+<0in,-.1in>*{\scriptstyle D_6}\restore &
{\bullet}\save[]+<0in,-.1in>*{\scriptstyle E_7}\restore
}} \\
&\vcenter{\xymatrix @!0 @M=0pt @R-.1in {
{\rlap{\text{(Type $E_8$)}}}&{}&{}&{}&{}&
{\circ}\save[]+<0in,-.1in>*{\scriptstyle D_4}\restore  \ar@{.}[r]&
{\circ}\save[]+<0in,-.1in>*{\scriptstyle D_6}\restore  \ar@{.}[r]&
{\bullet}\save[]+<0in,-.1in>*{\scriptstyle E_7}\restore  \ar@{-}[r]&
{\bullet}\save[]+<0in,-.1in>*{\scriptstyle E_8}\restore
}} \\
&\vcenter{\xymatrix  @!0 @M=0pt @R-.1in {
&{}&{}&{}&{}&{}&{\bullet}\save[]+<0in,-.1in>*{\scriptstyle A_{1}}\restore \ar@{-}'[dr] [ddrr] \\
&{}&{}&{}&{}&{}&{}&{\bullet}\save[]+<0in,-.1in>*{\scriptstyle C_3}\restore \\
{\rlap{\text{(Type $F_4$)}}}&{}&{}&{}&{}&{}&
{\circ}\save[]+<-.15in,-.1in>*{\scriptstyle B_2=C_2}\restore \ar@{.}[ur]
\ar@{.}[dr] &
{}& {\bullet}\save[]+<0in,-.1in>*{\scriptstyle F_4}\restore \\
&{}&{}&{}&{}&{}&{}&{\bullet}\save[]+<0in,-.1in>*{\scriptstyle B_3}\restore \\
&{}&{}&{}&{}&{}&{\bullet}\save[]+<0in,-.1in>*{\scriptstyle A_{1}}\restore \ar@{-}'[ur] [uurr]
}} \\
&\vcenter{\xymatrix @!0 @M=0pt @R-.1in {
&{}&{}&{}&{}&{\bullet}\save[]+<0in,-.1in>*{\scriptstyle A_{1}}\restore \ar@{-}[dr]\\
{\rlap{\text{(Type $G_2$)}}}&{}&{}&{}&{}&{}&{\bullet}\save[]+<0in,-.1in>*{\scriptstyle G_2}\restore \\
&{}&{}&{}&{}&{\bullet}\save[]+<0in,-.1in>*{\scriptstyle A_{1}}\restore \ar@{-}[ur]
}}
\end{align*}
\caption{Hasse Diagrams of $\Ftilde_C$ ($\bullet=$ node of $\F_C$; $\circ=$
 node of $\Ftilde_C\setminus \F_C$).  Each node $\psi$ is labeled by the
 type of the root system generated by $\psi$; for simplicity, $A_1$ nodes
 in $\Ftilde_C\setminus \F_C$ are omitted.}
\label{figHasseDiagrams}
\end{figure}
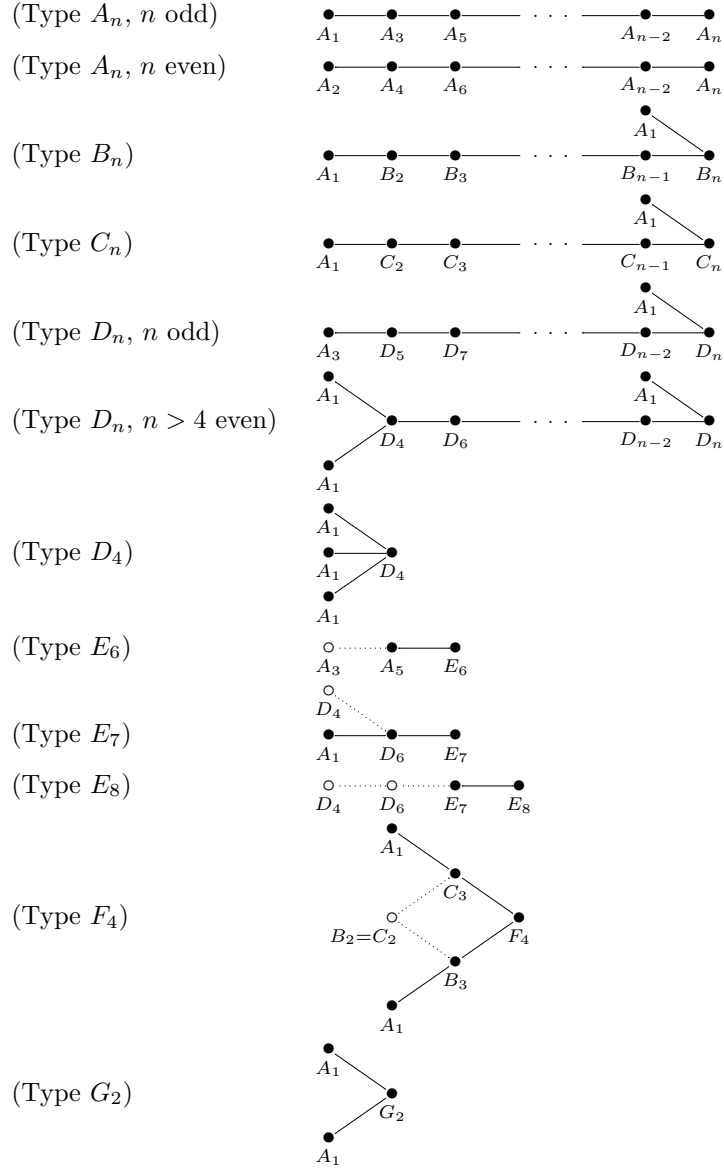
\end{proof}

\begin{defn}
\label{defnB}
Given a  Satake compactification  $\Dstar_\pi$, 
let $\B$ denote the family of subsets $\kap(\e_{\CC/\RR}(\theta))\subseteq
\lsb\CC\D$, where $\theta\subseteq \lsb\RR\D$ is a nonempty connected and
$\d$-connected subset. Let $\B^*\subseteq \B$ exclude those subsets
$\kap(\e_{\CC/\RR}(\theta))$ with cardinality $1$ which are not components
of $\lsb\CC\D$.  For every component $C$ of $\lsb\CC\D$, set $\B_{C\cup
c^*C} = \{\psi\in \B\mid \psi\subseteq C\cup c^*C\}$ and $\B^*_{C\cup c^*C}
= \B_{C\cup c^*C}\cap\B^*$.
\end{defn}

An element of $\B$ is simply the set of simple roots for the automorphism
group of a nontrivial standard irreducible boundary component.  The
correspondence $\theta \leftrightarrow \kap(\e_{\CC/\RR}(\theta))$ above is
an isomorphism of partially ordered sets for the inclusion ordering.  Note
that in the situation to be considered below, $C\cup c^* C= C$.

A semisimple algebraic $\RR$-group is \emph{equal-rank} if $\CCrank G =
\rank K$.  A symmetric space $D$ is \emph{equal-rank} if $D= G/K$ for an
equal-rank group $G$.  A \emph{real equal-rank Satake compactification} is
a Satake compactification $\Dstar_\pi$ for which all real
boundary components $D_{P,h}$ are equal-rank.  In the case of a real
equal-rank Satake compactification we wish to relate $\B$ (which depends on
the $\RR$-structure of $G$ and the given Satake compactification) with $\F$
(which only depends on $G$ as a $\CC$-group).  We begin with a basic lemma.

\begin{lem}[Borel and Casselman]
\label{lemEqualRank}
A semisimple $\RR$-group $G$ is equal-rank if and only if $c^*=\iota$,
where $c\in \Gal(\CC/\RR)$ is complex conjugation.
\end{lem}
\begin{proof}
See \cite[\S1.2(1) and Corollary~ 1.6(b)]{refnBorelCasselman}.
\end{proof}

\begin{cor}
\label{corConjugationCommutes}
If a semisimple $\QQ$-group $G$ is equal-rank, then $c^*$ commutes with
$g^*$ for all $g\in \Gal(\CC/\QQ)$.
\end{cor}

\begin{lem}
\label{lemBstarInFstar}
Let $\Dstar_\pi$ be a real equal-rank Satake compactification.  Then $\B
\subseteq \Ftilde$, every component of $\D^0_{\CC/\RR}$ belongs to
$\Ftilde$, and every component $C$ of $\lsb\CC\D$ is $c^*$-stable.
\end{lem}
\begin{proof}
Let $c^*_\psi$ denote the $*$-action of $c$ for the subroot system with
simple roots $\psi=\kap(\e_{\CC/\RR}(\theta))$; one may check that
$c^*_\psi = c^*|_\psi$.  By our equal-rank assumption and Lemma
~\ref{lemEqualRank}, $c^*_\psi = \iota_\psi$ and $c^* =\iota$.  It follows
that $\psi$ is $\iota$-stable and that $\iota|_\psi = \iota_\psi$.  It also
follows that $\psi$ is connected, since otherwise $c^*_\psi$ would
interchange the components and $\iota_\psi$ would preserve them.  The same
argument applies to a component $\psi$ of $\D^0_{\CC/\RR}$ since $c^*_\psi
= c^*|_\psi$ and $\psi$ is the index of an equal-rank (even
$\RR$-anisotropic) group.
\end{proof}

\begin{cor}
\label{corFFour}
Let $\Dstar_\pi$ be a real equal-rank Satake compactification.  If a
component $C$ of $\lsb\CC\D$ has type $F_4$, then $\RRrank C = 0$ or $1$.
\end{cor}
\begin{proof}
If $\RRrank C = 4$, then $\B_C\subseteq \Ftilde_C$ contains a chain of
length $4$.  By Figure ~\ref{figHasseDiagrams}, the chain must begin with
one of the roots in the $B_2$ subdiagram; this root must therefore belong
to $\d$.  But then $\B_C$ would also contain an element of type $A_2$ which
is excluded by the figure.  Thus $\RRrank C \le 3$.  However $F_4$ only has
real forms with $\RRrank= 0$, $1$, and $4$ \cite[\S5.9]{refnAraki},
\cite[Chapter X, Table V]{refnHelgason}.
\end{proof}

\begin{cor}
\label{corConnected}
Let $\Dstar_\pi$ be a real equal-rank Satake compactification.  If
$\theta\subseteq \lsb\RR C$ is a $\d$-connected subset of a connected
component $\lsb\RR C$ of $\lsb\RR\D$, then $\theta$ is connected.
Consequently $\kap(\e_{\CC/\RR}(\theta)) \subseteq \lsb\CC\D$ is connected
and \textup(if nonempty\textup) belongs to $\B$.
\end{cor}
\begin{proof}
Let $C = \e_{\CC/\RR}(\lsb\RR C)$ be the corresponding component of
$\lsb\CC\D$; it is connected by Lemma ~\ref{lemBstarInFstar}.  We will
prove that any two elements of $\lsb\RR\d\cap\lsb\RR C$ are connected by an
edge,
from which the corollary follows.  To see this, assume otherwise.  Then
there are clearly two distinct connected and $\d$-connected subsets of
$\lsb\RR C$ with
cardinality $2$.  This is impossible if the type of $C$ is not $F_4$, since
$\B^*_C\subseteq \Ftilde^*_C$ is totally ordered by Proposition
~\ref{propFTotallyOrdered}\itemref{itemNotFFour} and Lemma
~\ref{lemBstarInFstar}.  Type $F_4$ is excluded by Corollary ~\ref{corFFour}.
\end{proof}

We introduce the following convenient notational convention: for $i\ge 1$,
$\psi_i$ will always denote some element of $\B$ with $\RRrank \psi_i = i$;
let $\psi_0=\emptyset$.  The following observation will be used repeatedly.
\begin{lem}
\label{lemCrucial} 
Let $\Dstar_\pi$ be a real equal-rank Satake compactification.
For any $\psi_{i-1} < \psi_i$ we may
decompose
\begin{equation}
  \label{eqnJump}
  \psi_i\setminus \psi_{i-1}= \restr_{\CC/\RR}^{-1}(\al) \cup \eta
  \qquad\text{for some $\al \in \lsb\RR\D$ and
    $\eta\subseteq \D^0_{\CC/\RR}$.}
\end{equation}
The subset $\eta$ is a union of components of $\D^0_{\CC/\RR}$ and
\begin{equation}
  \label{eqnConnections}
  \psi_{i-1} \cup \{\b\} \text{ is }
  \begin{cases}
    \text{connected} & \text{for $\b\in \restr_{\CC/\RR}^{-1}(\al)$,}\\
    \text{disconnected}   & \text{for $\b\in \eta$ if $i>1$.}
  \end{cases}
\end{equation}
\end{lem}
\begin{proof}
This follows from Definition ~\ref{defnB} except for the assertion of
``connected'' in
\eqref{eqnConnections} as opposed to ``$\d$-connected''; for that use (the
proof of) Corollary ~\ref{corConnected}.
\end{proof}

\begin{prop}
\label{propBTotallyOrdered}
Let $\Dstar_\pi$ be a real equal-rank Satake compactification.
Then $\B\subseteq \F$.
\end{prop}

\begin{proof}
By Lemma ~\ref{lemBstarInFstar} it suffices to verify that any
$\psi=\psi_i\in \B_C$ satisfies conditions \itemref{itemFone} and
\itemref{itemFtwo} from Definition~\ref{defnF}.  If $\psi_i^+\setminus
\psi_i$ had two elements modulo $\iota$, there would exist two incomparable
elements $\psi_{i+1}$ and $\psi'_{i+1}$ of $\B^*_C$ strictly containing
$\psi_i$.  For $C$ not type $F_4$, this is impossible by Lemma
~\ref{lemBstarInFstar} and Proposition
~\ref{propFTotallyOrdered}\itemref{itemNotFFour} since $\B_C^*\subseteq
\Ftilde_C^*$ is totally ordered.    Type $F_4$ is excluded by Corollary
~\ref{corFFour}.  Condition 
\itemref{itemFone} follows.

As for condition \itemref{itemFtwo}, let $\psi_{i+1}\in \B_C$ contain
$\psi_i$; \itemref{itemFone} implies that $\psi_{i+1}\supseteq \psi_i^+$.
Then $\psi_{i+1}\setminus\psi_i^+ =\eta$ is a union of components of
$\D^0_{\CC/\RR}$ by Lemma ~\ref{lemCrucial}; these components belong to
$\Ftilde_C$ by Lemma ~\ref{lemBstarInFstar}.
\end{proof}

\begin{cor}
\label{corOneEnd}
Let $\Dstar_\pi$ be a real equal-rank Satake compactification.
Consider $\psi_{i-1} < \psi_i \in \B$ and let $\al\in\lsb\RR\D$ be as in
\eqref{eqnJump}.  If $\psi = \psi_{i-1}$ or if
$\psi\in \F$ is a component of $(\psi_i\setminus \psi_{i-1})\cap
\D^0_{\CC/\RR}$ \textup(the \emph{noncompact} case and the \emph{compact}
case respectively\textup), then $\psi^+\setminus \psi =
\restr_{\CC/\RR}^{-1}(\al)$.
\end{cor}
\begin{proof}
In the noncompact case, \eqref{eqnConnections} implies
$\restr_{\CC/\RR}^{-1}(\al) \subseteq \psi^+\setminus \psi$; the same
inclusion holds in the compact case since $\psi_i$ is connected.  However
$\psi\in\F$ by the proposition (or by hypothesis in the compact case).
Thus Definition~\ref{defnF}\itemref{itemFone} and Lemma~\ref{lemEqualRank}
imply that $\psi^+\setminus \psi$ is a single $c^*$-orbit.  The corollary
follows.
\end{proof}

\begin{thm}
\label{thmMainTheorem}
Let $G$ be an almost $\QQ$-simple semisimple group and let
$\Dstar_\pi$ be a real equal-rank Satake compactification.  If
$G$ has an $\RR$-simple factor $H$ with $\RRrank H=2$ and $\CC$-type
$B_n$, $C_n$, or $G_2$, assume that the Satake compactification associated
to $\pi|_{H}$ does not have a real boundary component of type $A_1$.
Then $\Dstar_\pi$ is geometrically rational.
\end{thm}
\begin{rem}
In particular, this implies Baily and Borel's result \cite{refnBailyBorel}
on geometric rationality for the natural compactification of a Hermitian
symmetric space except possibly when $G$ is the restriction of scalars of a
group with $\CC$-type $B_n$, $C_n$, or $G_2$.  But the Hermitian
condition excludes $G_2$ and implies for $B_n$ or $C_n$ that all simple
factors of $G$ have the same $\RR$-type, so $\d$ is Galois invariant and
Theorem ~\ref{thmDeltaRational} applies.
\end{rem}

\begin{proof}
For every component $C$ of $\lsb\CC\D$, define
\begin{equation*}
  \F_C^\circ = \F_C\setminus \left\{\psi\in\F_C \left|
  \begin{aligned}
    &\text{$\psi$ is type $A_1$ and is covered by}\\
    &\text{$C$ which is type $B_n$, $C_n$, or $G_2$}
  \end{aligned}
  \right.\right\}
\end{equation*}
and set $\F^\circ =\coprod_C \F^\circ_C$.  If $\psi\in \F_C\setminus
\F_C^\circ$ belonged to $\B$, this would imply that $\RRrank C = 2$.  Since
such a component is excluded by our hypotheses, Proposition
~\ref{propBTotallyOrdered} may be strengthened to $\B\subseteq \F^\circ$.
Proposition ~\ref{propFTotallyOrdered}\itemref{itemIncomparable} implies
that
\begin{equation}
\label{eqnIncomparable}
  \text{if $\psi\in\B_C$ and $\psi'\in\F^\circ_C$ are incomparable, then
  $\psi\cup\psi'$ is disconnected.}
\end{equation}
(If $C$ has type $F_4$, $\RRrank C = 1$ by Corollary
~\ref{corFFour} which implies $\B_C = \{C\}$ and the above
equation is vacuous.)  It follows that
\begin{equation}
\label{eqnBTotallyOrdered}
  \text{$\B_C$ is totally ordered.}
\end{equation}
For if $\psi$, $\psi'\in\B_C$ were incomparable, the union would be both
disconnected by \eqref{eqnIncomparable} and connected by Corollary
~\ref{corConnected}.

Define
\begin{align*}
&\K_{nc} =  \{\,\psi \mid \text{$\psi$ is a component of
  $\kap(\D^0_{\CC/\QQ})$}\,\}, \\
&\K_c = \left\{\,\psi\in\F^\circ\left|
   \begin{aligned}
     &\text{$\psi$ is a component of
     $(\psi_i\setminus\psi_{i-1})\cap\D^0_{\CC/\RR}$ for}\\ 
     &\qquad\text{some  $\psi_{i-1}<\psi_i\in\B$ with $\psi_i\not\subseteq
     \kap(\D^0_{\CC/\QQ})$}
   \end{aligned}
\right.\,\right\},
\end{align*}
and set $\K = \K_{nc}\coprod \K_c$; we call the elements of $\K_{nc}$
\emph{noncompact} and the elements of $\K_c$ \emph{compact}.  For a
component $C$ of $\lsb\CC\D$, we define $\K_C$, $\K_{C,nc}$, and $\K_{C,c}$
as usual.

Note that $\kap(\D^0_{\CC/\QQ}) = \kap(\e_{\CC/\RR}(\D^0_{\RR/\QQ})) =
\kap(\e_{\CC/\RR}(\kap(\D^0_{\RR/\QQ})))$; thus by Corollary
~\ref{corConnected}
\begin{equation}
\label{eqnKNonCompact}
\K_{nc} =  \{\,\kap(\D^0_{\CC/\QQ})\cap C \mid \d \cap
  \D^0_{\CC/\QQ} \cap C \neq \emptyset \,\} \subseteq \B
  \subseteq \F^\circ.
\end{equation}

We first show that if $\psi\in \K$ and $g\in \Gal(\CC/\QQ)$, then
\begin{equation}
\label{eqnKProperties}
  \text{$g^*\psi\in\F^\circ$, $g^*\psi\subseteq \D^0_{\CC/\QQ}$, and
  $g^*\psi$ is maximal among such sets.}
\end{equation}
It suffices to prove this for $g^*$ the identity since the definition of
$\F^\circ$ depends only the $\CC$-root system and since $\D^0_{\CC/\QQ}$ is
$\Gal(\CC/\QQ)$-invariant.  Then the first assertion is part of the
definition if $\psi$ is compact and follows from \eqref{eqnKNonCompact} if
$\psi$ is noncompact.  The second assertion, that $\psi\subseteq
\D^0_{\CC/\QQ}$, is clear.  For the final assertion, assume $\psi\in\K_C$
and suppose that there exists $\psi'\in\F^\circ$ such that $\psi\subset
\psi' \subseteq \D_{\CC/\QQ}^0 \cap C$.  If $\psi$ is noncompact then it is
$\d$-connected and hence $\psi'$ (being connected) must also be
$\d$-connected; this contradicts the fact that
$\psi=\kap(\D^0_{\CC/\QQ})\cap C$ is the largest $\d$-connected subset of
$\D_{\CC/\QQ}^0\cap C$.  If instead $\psi\subseteq
\psi_i\setminus\psi_{i-1}$ is compact, let $\al\in\lsb\RR\D$ be as in
\eqref{eqnJump}.  Since $\psi^+\setminus \psi = \restr_{\CC/\RR}^{-1}(\al)$
by Corollary \ref{corOneEnd}, any connected set
strictly containing $\psi$, such as $\psi'$, must contain an
element of $\restr_{\CC/\RR}^{-1}(\al)$.  Thus (a)
$\psi_i\setminus\psi_{i-1}\subseteq \D^0_{\CC/\QQ}$ by \eqref{eqnJump} and
(b) $\psi_{i-1}\cup
\psi'$ is connected by \eqref{eqnConnections}.  If $i=1$, (a)
implies $\psi_1\subseteq \kap(\D^0_{\CC/\QQ})$ which contradicts the
definition of $\K_c$.  If
$i>1$, the same argument shows that $\psi_{i-1} \nsubseteq \psi'$ and hence
$\psi_{i-1}$ and $\psi'$ are incomparable; this contradicts
\eqref{eqnIncomparable} and (b).

We now prove  that
\begin{equation}
\label{eqnKGaloisInvariant}
\text{$\K$ is $\Gal(\CC/\QQ)$-invariant,}
\end{equation}
that is, given $\psi\in \K_C$ and $g\in \Gal(\CC/\QQ)$ we will show that
$g^*\psi\in \K$.  Since we can assume $G$ is not $\QQ$-anisotropic,
$g^*\psi\neq g^*C$.  Let $i$ be minimal such that $g^*\psi < \psi_i$ for
some $\psi_i\in\B$ and consider some $\psi_{i-1} < \psi_i$.  We know that
$\psi_i\not \subseteq \kap(\D^0_{\CC/\QQ})$, else the maximality in
\eqref{eqnKProperties} would be contradicted.  There are three cases.

\emph{Case} 1: $i=1$.  Let $\al$ and $\eta$ be as in  Lemma
~\ref{lemCrucial}.  We have $\restr_{\CC/\RR}^{-1}(\al) \cap
\D^0_{\CC/\QQ} = \emptyset$ else $\psi_1\subseteq \kap(\D^0_{\CC/\QQ})$.
Thus $g^*\psi\in\F^\circ$ is contained in a component of $\eta =
\psi_1\cap \D^0_{\CC/\RR}$; since Proposition ~\ref{propFTotallyOrdered}
and Lemmas ~\ref{lemBstarInFstar} and \ref{lemCrucial} imply such a
component is in $\F^\circ$, and hence $\K_c$, the maximality in
\eqref{eqnKProperties} implies $g^*\psi$ equals this
component.

\emph{Case} 2: $i>1$ and $\psi_{i-1}\le g^*\psi$.  In this case,
$\psi_{i-1}\subseteq g^*\psi\subseteq \D_{\CC/\QQ}^0$ where $\psi_{i-1}$ is
nonempty and $\d$-connected and $g^*\psi$ is connected.  This implies
$g^*\psi$ is $\d$-connected and hence equals a component of
$\kap(\D_{\CC/\QQ}^0)$ by the maximality in \eqref{eqnKProperties}.  So
$g^*\psi\in\K_{nc}$.

\emph{Case} 3: $i>1$ and $\psi_{i-1}\not\le g^*\psi$.  Since $\psi_{i-1}
\cup g^*\psi$ is disconnected by \eqref{eqnIncomparable}, $g^*\psi$ is
contained in a connected component of
$(\psi_i\setminus\psi_{i-1})\cap\D^0_{\CC/\RR}$ by \eqref{eqnConnections};
as in Case 1, this component must be in $\K_c$ and $g^*\psi$ equals it.

This finishes the proof that $\K$ is $\Gal(\CC/\QQ)$-invariant.

Let $\widehat\K = \Gal(\CC/\QQ)\cdot \K_{nc}$. We claim that if $\psi \in
\widehat\K_C$ is compact, say $\psi \subseteq
(\psi_i\setminus\psi_{i-1})\cap\D^0_{\CC/\RR}$, then $\psi_{i-1} =
\kap(\D^0_{\CC/\QQ})\cap C$.  In the case $i=1$, the claim asserts that
$\kap(\D^0_{\CC/\QQ})\cap C=\emptyset$; this holds since otherwise
$\kap(\D^0_{\CC/\QQ})\cap C\supseteq \psi_1$ (by
\eqref{eqnBTotallyOrdered}) which contradicts the definition of $\K_c$. As
for the case $i>1$, let $g\in\Gal(\CC/\QQ)$ be such that $g^*\psi$ is
noncompact.  Then $g^*\psi\in\B$ and $g^*\psi_{i-1}$ are incomparable; let
$\tilde \psi\in\B$ contain $g^*\psi$ and be maximal such that $\tilde\psi$
and $g^*\psi_{i-1}$ are incomparable.  It follows from
\eqref{eqnIncomparable} and Lemma ~\ref{lemCrucial} that
$g^*\psi_{i-1}\subseteq \D^0_{\CC/\RR}$ and hence $\psi_{i-1} \subseteq
\D^0_{\CC/\QQ}$.
Thus $\psi_{i-1}\subseteq \kap(\D^0_{\CC/\QQ})\cap C\in \B_C$, while $\psi_i
\not\subseteq \kap(\D^0_{\CC/\QQ})$.  The claim now follows from
\eqref{eqnBTotallyOrdered}.

Assume now that $\widehat\K_C\neq\emptyset$ for one and hence all
components $C$.  For a fixed $C$, it follows from the above claim and
Corollary ~\ref{corOneEnd} that $\psi^+\setminus \psi$ is independent of
the choice of $\psi\in\widehat \K_C$ and, if there does not exist a
noncompact element in $\widehat\K_C$, then $\d\cap C= \psi^+\setminus\psi$.
It follows that
\begin{equation}
\label{eqnOmegaCondition}
  \b\in\o(\D^0_{\CC/\QQ})\cap C \quad\Longleftrightarrow\quad \text{
    $\b\notin \psi^+\setminus\psi$ for some (and hence any)
    $\psi\in\widehat\K_C$.}
\end{equation}
For in the case that there exists $\psi \in \widehat\K_C$ noncompact,
equation \eqref{eqnOmegaCondition} holds since $\psi =
\kap(\D^0_{\CC/\QQ})\cap C$, while in the case that every $\psi \in
\widehat\K_C$ is compact, $\o(\D^0_{\CC/\QQ})\cap C$ is simply the
complement of $\d\cap C = \psi^+\setminus\psi$ in $C$.  Since the right
hand side of \eqref{eqnOmegaCondition} is $\Gal(\CC/\QQ)$-invariant,
condition \itemref{itemGRone} of Theorem ~\ref{thmCasselmanCriterion}
holds.  Also if $g\in\Gal(\CC/\QQ)$ and $\psi\in \K$ is a component of
$\kap(\D^0_{\CC/\QQ})$, then $g^*\psi\in \K$ and thus either
$g^*\psi\subseteq \kap(\D^0_{\CC/\QQ})$ or $g^*\psi \subseteq
\D^0_{\CC/\RR}$; this verifies condition \itemref{itemGRtwo} of Theorem
~\ref{thmCasselmanCriterion}.

The remaining case is where $\widehat\K_C=\emptyset$ for all $C$.  Then
$\bigcup_{g\in\Gal(\CC/\QQ)} g^*\B\cap C$ is totally ordered for any $C$;
this may be proved similarly to \eqref{eqnKGaloisInvariant} but more
easily.  Let $\psi_1\in \B_C$ and compare $g^*\psi_1$ with $\psi'_1 \in
\B_{g^*C}$.  Say $g^*\psi_1 \le \psi'_1$.  We must have $g^*(\d\cap C)
\subseteq \d\cap g^*C$, since otherwise $\kap(\D^0_{\CC/\QQ})$ would not be
empty.  Since $g^*$ commutes with $c^*$ (Corollary
~\ref{corConjugationCommutes}) and $\d\cap C$ has only one element modulo
$c^*$ (by equation \eqref{eqnBTotallyOrdered}), we have equality.  Thus
$\d$ is Galois invariant and the theorem follows from Theorem
~\ref{thmDeltaRational}.
\end{proof}

\section{Exceptional cases}
\label{sectExceptionalCases}
It remains to handle the cases excluded from Theorem ~\ref{thmMainTheorem}.
In contrast to the situation of that theorem, where geometric rationality
was automatic with no rationality assumption on $(\pi,V)$, these
exceptional cases are 
only geometrically rational under a certain condition which depends on the
$\QQrank$.

\begin{thm}
\label{thmSpecialCases}
Let $G$ be an almost $\QQ$-simple semisimple group which is the restriction
of scalars of a group $G'$ with $\CC$-type $B_n$, $C_n$, or $G_2$.  Let
$\Dstar_\pi$ be a real equal-rank Satake compactification.
Assume $G$ has an $\RR$-simple factor $H$ with $\RRrank H=2$ for which the
Satake compactification associated to $\pi|_{H}$ has a real boundary
component of type $A_1$.
\begin{enumerate}
\item \label{itemQrankTwo} In the case $\QQrank G=2$,
$\Dstar_\pi$ is geometrically rational if and only if $\d$ is
Galois invariant.
\item \label{itemQrankOne} In the case $\QQrank G=1$,
$\Dstar_\pi$ is geometrically rational if and only if $\d\cap
\D^0_{\CC/\QQ}$ is empty or meets every component of $\lsb\CC\D$ with
$\RRrank \ge 2$.
\end{enumerate}
\end{thm}

\begin{proof}
Let $C$ be the component of $\DCC$ corresponding to $H$ and let $\al_1\in
\d\cap C$ correspond to the real boundary component of type $A_1$.  By
Lemma ~\ref{lemCrucial} and our hypotheses, $\al_1$ is at an end of $C$ and
(if we denote its unique neighbor by $\al_2$) we have $C\setminus \{\al_1,
\al_2\}= \D^0_{\CC/\RR}\cap C\subseteq \D^0_{\CC/\QQ}\cap C$.  If $H$ has
$\CC$-type $C_n$, the classification of semisimple Lie algebras over $\RR$
\cite{refnAraki}, \cite[Chapter X, Table VI]{refnHelgason}
shows that the two adjacent non-$\RR$-anisotropic roots imply $H$ is
$\RR$-split; thus $C_n$ only occurs for $n=2$ and we can absorb this case
into that of $B_n$.  Classification theory over $\RR$ also shows that if
$n>2$, the simple root $\al_1$ is long.

For \itemref{itemQrankTwo} we only need to prove that $\d$ is Galois
invariant under the assumption that \itemref{itemGRone} and
\itemref{itemGRtwo} of Theorem ~\ref{thmCasselmanCriterion} hold; the
opposite direction was already proved in Theorem ~\ref{thmDeltaRational}.
Thus assume that $\QQrank G=2$.  Then for all $g\in\Gal(\CC/\QQ)$ we have
$\kap(\D^0_{\CC/\QQ})\cap g^*C \subseteq \D^0_{\CC/\QQ}\cap g^*C =
g^*(C\setminus \{\al_1, \al_2\}) \subseteq \o(\D^0_{\CC/\QQ})\cap g^*C$.
Since $g^*\al_1$ is orthogonal to $g^*(C\setminus \{\al_1, \al_2\})$, it
follows that $g^*\al_1\notin\o(\D^0_{\CC/\QQ})\cap g^*C$ if and only if
$g^*\al_1 \in \d\cap g^*C$.  But the first condition is independent of $g$
by condition \itemref{itemGRone} of Theorem ~\ref{thmCasselmanCriterion}
while the second condition holds for $g=e$, so we must have $g^*\al_1 \in
\d\cap g^*C$ for all $g\in\Gal(\CC/\QQ)$.
If $\#(\d \cap g_0^*C) > 1$ for some $g_0\in \Gal(\CC/\QQ)$ then $\RRrank
g_0^* C = 2$ by Propositions
~\ref{propFTotallyOrdered}\itemref{itemIncomparable} and
\ref{propBTotallyOrdered}, in which case $g_0^*(C\setminus \{\al_1,
\al_2\})= \D^0_{\CC/\RR}\cap g_0^*C$; thus the additional element of $\d
\cap g_0^*C$ must be $g_0^*\al_2$.  In this case however
$\kap(\D^0_{\CC/\QQ}) = \emptyset$ so that for all $g\in \Gal(\CC/\QQ)$,
$g^*\al_2\notin\o(\D^0_{\CC/\QQ})\cap g^*C$ if and only if $g^*\al_2 \in
\d\cap g^*C$.  Another application of condition \itemref{itemGRone} from
Theorem ~\ref{thmCasselmanCriterion} concludes the proof that $\d$ is
Galois invariant in this case.

For \itemref{itemQrankOne} assume that $\QQrank G=1$.  Write $G = R_{k/\QQ}
G'$ where $k$ is a finite extension of $\QQ$ (totally real by Lemma
~\ref{lemBstarInFstar}) and $G'$ is an almost $k$-simple group with
$k\fieldrank=1$ and $\CC$-type $B_n$ or $G_2$.  The classification of
semisimple groups over $\QQ$ \cite{refnTits} shows the case $G_2$ cannot
occur and that in the case $B_n$ the unique simple root which is not
$k$-anisotropic is the long root at the end of the Dynkin diagram.
Since $\D^0_{\CC/\QQ}\cap C$ is thus connected for every component $C$ of
$\lsb\CC\D$, it follows that
\begin{equation}
\text{$\kap(\D^0_{\CC/\QQ})\cap C$ is nonempty}
    \qquad\Longleftrightarrow\qquad\text{$\d\cap \D^0_{\CC/\QQ}$ meets
    $C$.}
\end{equation}
If this condition holds, then $\kap(\D^0_{\CC/\QQ})\cap
C=\D^0_{\CC/\QQ}\cap C$, while if the condition does not hold, then $\d\cap
C = C\setminus (\D^0_{\CC/\QQ}\cap C)$.  In either case,
$\o(\D^0_{\CC/\QQ})\cap C=\D^0_{\CC/\QQ}\cap C$.  Furthermore
\begin{equation}
\D^0_{\CC/\QQ}\cap C \not\subseteq \D^0_{\CC/\RR}
\qquad\Longleftrightarrow\qquad
\RRrank C\ge 2.
\end{equation}
Part \itemref{itemQrankOne} now easily follows from Theorem
~\ref{thmCasselmanCriterion}.
\end{proof}

\begin{rems}
\label{remExceptionalCases}
Assume $G$ satisfies the hypotheses of Theorem ~\ref{thmSpecialCases}.
\begin{enumerate}
\item \label{itemQrankTwoRemark} In the $\QQrank 2$ case, if
  $\Dstar_\pi$ is geometrically rational then all $\RR$-simple
  factors of $G$ have $\RRrank = 2$.  For the proof shows that every
  component $g^*C$ has an $A_1$ boundary component corresponding to
  $g^*\al_1\in \d\cap g^*C$ and that if $n>2$ then $\al_1$ is long.
  However if $\RRrank g^*C>2$ then the entry for $B_n$ in Figure
  ~\ref{figHasseDiagrams} implies that $\d\cap g^*C$ is a singleton short
  root.
\item \label{itemQrankOneRrankBigRemark} In the $\QQrank 1$ case, $\d\cap
  \D^0_{\CC/\QQ}$ automatically meets any component with $\RRrank > 2$;
  this follows from the proof and the entry for $B_n$ in Figure
  ~\ref{figHasseDiagrams}.  Thus if there is any component with $\RRrank >
  2$ , the condition for geometric rationality is that $\d\cap
  \D^0_{\CC/\QQ}$ meets every component with $\RRrank = 2$.
\item \label{itemQrankOneRrankOneRemark} Also in the $\QQrank 1$ case, if
  $\CCrank G' > 2$ then $G$ has a component with $\RRrank 1$.  For by the
  proof above, $G'$ (up to strict $k$-isogeny) is the special orthogonal
  group of a quadratic form $q$ in at least $7$ variables with $k$-index
  $1$.  Thus we can decompose $q$ into a hyperbolic plane and a
  $k$-anisotropic form $q'$ in at least $5$ variables.  Since every
  quadratic form in at least $5$ variables is isotropic at every
  non-archimedean place \cite[VI.2.12]{refnLamQuadraticForms}, $q'$ must be
  anisotropic at some real place by the Hasse-Minkowski principle
  \cite[VI.3.5]{refnLamQuadraticForms}.
\end{enumerate}
\end{rems}

\begin{examples}
Let $k$ be a real quadratic extension of $\QQ$.  In the Satake diagrams
below, the roots in $\D^0_{\CC/\RR}$ are colored black, while roots in
$\D^0_{\CC/\QQ}$ are enclosed in a dotted box.  Nodes placed vertically
above each other form a $\Gal(\CC/\QQ)$-orbit.  The elements of $\d$ are so
labeled.

\begin{enumerate}
\item Let $q$ be a $k$-anisotropic form in three variables,
  $\RR$-anisotropic at one real place and of signature $(2,1)$ at the
  other.  Let $G'$ be the orthogonal group of $h\oplus h \oplus q$, where
  $h$ is hyperbolic space, and let $G=R_{k/\QQ}G'$.  Then both of the
  Satake compactifications
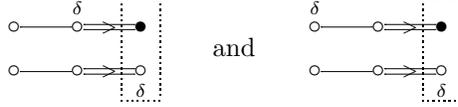
\begin{equation*}
\vcenter{\xymatrix @!0 @M=0pt @R-.1in {
{\circ} \ar@{-}[r] &
{\circ}\save[]+<0in,.1in>*{\scriptstyle \d}\restore \ar@{=}[r] \ar@{}[r]|{>}&
{\bullet}\save[]+<-.1in,.15in>="a"\restore\\
{\circ} \ar@{-}[r] &
{\circ} \ar@{=}[r] \ar@{}[r]|{>}&
{\circ}\save[]+<0in,-.1in>*{\scriptstyle \d}\restore
\save[]+<.1in,-.15in>="b"\restore
\save"a"."b"*[F.]\frm{}\restore
}}
\qquad\text{and}\qquad
\vcenter{\xymatrix @!0 @M=0pt @R-.1in {
{\circ}\save[]+<0in,.1in>*{\scriptstyle \d}\restore \ar@{-}[r] &
{\circ} \ar@{=}[r] \ar@{}[r]|{>}&
{\bullet}\save[]+<-.1in,.15in>="a"\restore\\
{\circ} \ar@{-}[r] &
{\circ} \ar@{=}[r] \ar@{}[r]|{>}&
{\circ}\save[]+<0in,-.1in>*{\scriptstyle \d}\restore
\save[]+<.1in,-.15in>="b"\restore
\save"a"."b"*[F.]\frm{}\restore
}}
\end{equation*}
are real equal-rank.  The first is geometrically rational (Theorem
~\ref{thmMainTheorem} applies) and the second is not (Theorem
~\ref{thmSpecialCases}\itemref{itemQrankTwo} applies).  More specifically,
condition \itemref{itemGRone} from Theorem ~\ref{thmCasselmanCriterion}
fails in the second example; this example also illustrates Remark
~\ref{remExceptionalCases}\itemref{itemQrankTwoRemark}.

If $q$ is instead nonsingular of dimension $1$, then of three real
equal-rank Satake compactifications
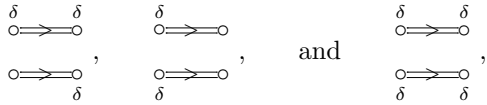
\begin{equation*}
\vcenter{\xymatrix @!0 @M=0pt @R-.1in {
{\circ}\save[]+<0in,.1in>*{\scriptstyle \d}\restore \ar@{=}[r]  \ar@{}[r]|{>}&
{\circ}\save[]+<0in,.1in>*{\scriptstyle \d}\restore\\
{\circ} \ar@{=}[r]  \ar@{}[r]|{>}&
{\circ}\save[]+<0in,-.1in>*{\scriptstyle \d}\restore 
}}
\text{ ,}\qquad
\vcenter{\xymatrix @!0 @M=0pt @R-.1in {
{\circ}\save[]+<0in,.1in>*{\scriptstyle \d}\restore \ar@{=}[r]  \ar@{}[r]|{>}&
{\circ} \\
{\circ}\save[]+<0in,-.1in>*{\scriptstyle \d}\restore \ar@{=}[r] \ar@{}[r]|{>}&
{\circ}
}}
\text{ ,}\qquad\text{and}\qquad
\vcenter{\xymatrix @!0 @M=0pt @R-.1in {
{\circ}\save[]+<0in,.1in>*{\scriptstyle \d}\restore \ar@{=}[r]  \ar@{}[r]|{>}&
{\circ}\save[]+<0in,.1in>*{\scriptstyle \d}\restore \\
{\circ}\save[]+<0in,-.1in>*{\scriptstyle \d}\restore \ar@{=}[r]  \ar@{}[r]|{>}&
{\circ}\save[]+<0in,-.1in>*{\scriptstyle \d}\restore
}}\text{ ,}
\end{equation*}
Theorem ~\ref{thmSpecialCases}\itemref{itemQrankTwo} shows that the first
is not geometrically rational and the last two are.  Similar examples can
be constructed with $q$ totally $\RR$-anisotropic in any odd number of
variables or where $G'$ is the split form of $G_2$.

\item
Let $G'$ be the orthogonal group of $h\oplus q$, where $q$ is a
$k$-anisotropic form in three variables which has signature $(2,1)$ at both
real places, and let $G=R_{k/\QQ}G'$.  Among the real equal-rank Satake
compactifications
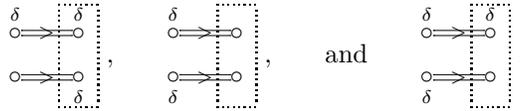
\begin{equation*}
\vcenter{\xymatrix @!0 @M=0pt @R-.1in {
{\circ}\save[]+<0in,.1in>*{\scriptstyle \d}\restore \ar@{=}[r] \ar@{}[r]|{>}&
{\circ}\save[]+<0in,.1in>*{\scriptstyle \d}\restore\save[]+<-.1in,.15in>="a"\restore\\
{\circ}\ar@{=}[r] \ar@{}[r]|{>}&
{\circ}
\save[]+<0in,-.1in>*{\scriptstyle \d}\restore 
\save[]+<.1in,-.15in>="b"\restore
\save"a"."b"*[F.]\frm{}\restore
}}
\text{ ,}\qquad
\vcenter{\xymatrix @!0 @M=0pt @R-.1in {
{\circ}\save[]+<0in,.1in>*{\scriptstyle \d}\restore  \ar@{=}[r] \ar@{}[r]|{>}&
{\circ}\save[]+<-.1in,.15in>="a"\restore\\
{\circ}\save[]+<0in,-.1in>*{\scriptstyle \d}\restore  \ar@{=}[r] \ar@{}[r]|{>}&
{\circ}
\save[]+<.1in,-.15in>="b"\restore
\save"a"."b"*[F.]\frm{}\restore
}}
\text{ ,}\qquad\text{and}\qquad
\vcenter{\xymatrix @!0 @M=0pt @R-.1in {
{\circ}\save[]+<0in,.1in>*{\scriptstyle \d}\restore \ar@{=}[r] \ar@{}[r]|{>}&
{\circ}\save[]+<0in,.1in>*{\scriptstyle \d}\restore 
\save[]+<-.1in,.15in>="a"\restore\\
{\circ}\save[]+<0in,-.1in>*{\scriptstyle \d}\restore 
\ar@{=}[r] \ar@{}[r]|{>}&
{\circ}
\save[]+<.1in,-.15in>="b"\restore
\save"a"."b"*[F.]\frm{}\restore
}}\text{ ,}
\end{equation*}
the first two are geometrically rational (even though in the first $\d$ is
not Galois invariant) and the last is not; this illustrates Theorem
~\ref{thmSpecialCases}\itemref{itemQrankOne}.

\item
Let $G'$ be the orthogonal group of $h\oplus q$, where $q$ is a
$k$-anisotropic form in five variables which has signature $(4,1)$ at one
real place and is $\RR$-anisotropic at the other (see Remark
\ref{remExceptionalCases}\itemref{itemQrankOneRrankOneRemark}) and let
$G=R_{k/\QQ}G'$.  There are two real equal-rank Satake compactifications
covered by Theorem ~\ref{thmSpecialCases}\itemref{itemQrankOne},
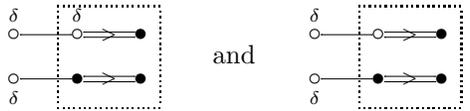
\begin{equation*}
\vcenter{\xymatrix @!0 @M=0pt @R-.1in {
{\circ}\save[]+<0in,.1in>*{\scriptstyle \d}\restore \ar@{-}[r] &
{\circ}\save[]+<0in,.1in>*{\scriptstyle \d}\restore
\save[]+<-.1in,.15in>="a"\restore \ar@{=}[r] \ar@{}[r]|{>}&
{\bullet}\\
{\circ}\save[]+<0in,-.1in>*{\scriptstyle \d}\restore \ar@{-}[r] &
{\bullet} \ar@{=}[r] \ar@{}[r]|{>}&
{\bullet}
\save[]+<.1in,-.15in>="b"\restore
\save"a"."b"*[F.]\frm{}\restore
}}
\qquad\text{and}\qquad
\vcenter{\xymatrix @!0 @M=0pt @R-.1in {
{\circ}\save[]+<0in,.1in>*{\scriptstyle \d}\restore \ar@{-}[r] &
{\circ}
\save[]+<-.1in,.15in>="a"\restore \ar@{=}[r] \ar@{}[r]|{>}&
{\bullet}\\
{\circ}\save[]+<0in,-.1in>*{\scriptstyle \d}\restore \ar@{-}[r] &
{\bullet} \ar@{=}[r] \ar@{}[r]|{>}&
{\bullet}
\save[]+<.1in,-.15in>="b"\restore
\save"a"."b"*[F.]\frm{}\restore
}}\text{ ,}
\end{equation*}
both of which are geometrically rational.  If, however, $k$ is a totally
real degree 3 extension of $\QQ$ and $q$ has signature $(4,1)$, $(5,0)$,
and $(3,2)$ at the three real places, the two real equal-rank Satake
compactifications covered by Theorem
~\ref{thmSpecialCases}\itemref{itemQrankOne} are
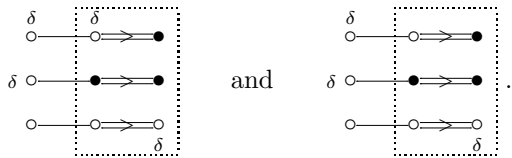
\begin{equation*}
\vcenter{\xymatrix @!0 @M=0pt @R-.1in {
{\circ}\save[]+<0in,.1in>*{\scriptstyle \d}\restore \ar@{-}[r] &
{\circ}\save[]+<0in,.1in>*{\scriptstyle \d}\restore
\save[]+<-.1in,.15in>="a"\restore \ar@{=}[r] \ar@{}[r]|{>}&
{\bullet}\\
{\circ}\save[]+<-.1in,0in>*{\scriptstyle \d}\restore \ar@{-}[r] &
{\bullet} \ar@{=}[r] \ar@{}[r]|{>}&
{\bullet}\\
{\circ} \ar@{-}[r] &
{\circ} \ar@{=}[r] \ar@{}[r]|{>}&
{\circ}\save[]+<0in,-.1in>*{\scriptstyle \d}\restore
\save[]+<.1in,-.15in>="b"\restore
\save"a"."b"*[F.]\frm{}\restore
}}
\qquad\text{and}\qquad
\vcenter{\xymatrix @!0 @M=0pt @R-.1in {
{\circ}\save[]+<0in,.1in>*{\scriptstyle \d}\restore \ar@{-}[r] &
{\circ}
\save[]+<-.1in,.15in>="a"\restore \ar@{=}[r] \ar@{}[r]|{>}&
{\bullet}\\
{\circ}\save[]+<-.1in,0in>*{\scriptstyle \d}\restore \ar@{-}[r] &
{\bullet} \ar@{=}[r] \ar@{}[r]|{>}&
{\bullet}\\
{\circ}\ar@{-}[r] &
{\circ} \ar@{=}[r] \ar@{}[r]|{>}&
{\circ}\save[]+<0in,-.1in>*{\scriptstyle \d}\restore 
\save[]+<.1in,-.15in>="b"\restore
\save"a"."b"*[F.]\frm{}\restore
}}\text{ .}
\end{equation*}
Here the first is geometrically rational and the second is not; this
illustrates Remark
~\ref{remExceptionalCases}\itemref{itemQrankOneRrankBigRemark}.
\end{enumerate}

\end{examples}


\bibliographystyle{amsplain}
\bibliography{greq}
\end{document}